# GEOMETRIC ISOMORPHISM AND MINIMUM ABERRATION FOR FACTORIAL DESIGNS WITH QUANTITATIVE FACTORS

By Shao-Wei Cheng[1] and Kenny Q. Ye[2]

*Academia Sinica and Albert Einstein College of Medicine*

Factorial designs have broad applications in agricultural, engineering and scientific studies. In constructing and studying properties of factorial designs, traditional design theory treats all factors as nominal. However, this is not appropriate for experiments that involve quantitative factors. For designs with quantitative factors, level permutation of one or more factors in a design matrix could result in different geometric structures, and, thus, different design properties. In this paper indicator functions are introduced to represent factorial designs. A polynomial form of indicator functions is used to characterize the geometric structure of those designs. *Geometric isomorphism* is defined for classifying designs with quantitative factors. Based on indicator functions, a new aberration criteria is proposed and some minimum aberration designs are presented.

**1. Introduction.** Factorial designs are commonly used in most industrial and scientific studies. In such a study, a number of fixed levels (settings) are selected for each factor (variable), and then some level combinations are chosen to be the runs in an experiment. A factor can be either nominal or quantitative. For nominal factors, there is no ordering among levels. The interest of analysis of an experiment with nominal factors is to understand if there exist differences in treatment means and if they exist, which treatment means differ. Analysis such as ANOVA or various multiple comparison testing procedures is often used for treatment comparison. In many studies, especially in response surface exploration, factors are often quantitative and there exists an order among levels. For an experiment with quantitative factors, the objective is usually achieved through fitting a (polynomial) model

Received August 2002; revised July 2003.
[1]Supported by the National Science Council of Taiwan, ROC.
[2]Supported by NSF Grant DMS-03-06306.

TABLE 1
*Combinatorially isomorphic designs with different geometric structures*

| $A$ | $B$ | $C$ | $A$ | $B$ | $C$ |
|---|---|---|---|---|---|
| 0 | 0 | 0 | 0 | 0 | 0 |
| 0 | 1 | 2 | 0 | 1 | 1 |
| 0 | 2 | 1 | 0 | 2 | 2 |
| 1 | 0 | 2 | 1 | 0 | 1 |
| 1 | 1 | 1 | 1 | 1 | 2 |
| 1 | 2 | 0 | 1 | 2 | 0 |
| 2 | 0 | 1 | 2 | 0 | 2 |
| 2 | 1 | 0 | 2 | 1 | 0 |
| 2 | 2 | 2 | 2 | 2 | 1 |

that can "well" describe the relationship between the response and the factors. The distinction in the analysis objective and strategy for these two types of experiments requires different selection criteria and classification methods.

For designs with nominal factors, the design properties should be invariant to level permutation within one or more of its factors. However, for quantitative factors, Cheng and Wu (2001) observed that level permutation of $3^{4-1}$ designs could result in changes in model efficiency when a polynomial model is fitted, which is referred to as "model nonisomorphism." Independently, Ye (1999) also observed that level permutation could alter the aliasing structure of designs when linear-quadratic decomposition [see Wu and Hamada (2000),

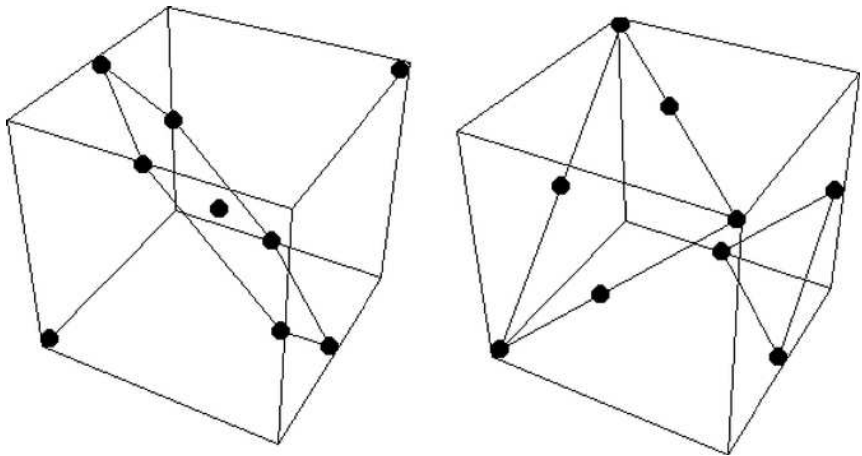

Fig. 1. *Combinatorially isomorphic but geometrically nonisomorphic designs.*

Section 5.6] is used. As shown in the following example, such "model nonisomorphism" is, indeed, the result of different geometric structures induced by permuting levels of factors. Consider the two $3^{3-1}$ designs in Table 1. In the table each design is written as a *design matrix* in which each column represents a factor and each row represents an experimental run. These designs are *combinatorially isomorphic* since one is obtained by applying the permutation $\{0,1,2\} \to \{0,2,1\}$ on the third column of the other. However, if we treat these levels as quantitative, their geometric structures are apparently different as shown in Figure 1. The difference in geometric structure also reflects on the model efficiency. For example, when a model that contains all linear main effects and three linear-by-linear interactions is considered, the design on the left-hand side has higher $D$-efficiency than the one on the right-hand side.

The conventional mathematical tools used for factorial designs, such as group theory and coding theory, treat all factors as nominal. Therefore, they do not differentiate geometric structures resulting from level permutations and fail to study the design properties associated with its geometric structure.

A new approach for characterizing designs with quantitative factors is developed in this paper. When all $k$ factors in a factorial design are quantitative, it can be viewed as a collection of points in $\mathbb{R}^k$. This collection of points is represented by an indicator function, which will be defined in Section 2. The indicator function can be written in a polynomial form which reveals the design's properties and characterizes its geometric structure. Thus, classification and design criteria are developed based on the indicator functions. This approach is motivated by Pistone and Wynn (1996), which first used polynomial systems to describe designs and studied their properties using algebraic geometry methods. In this paper properties of designs with quantitative factors are studied. Section 2 introduces the indicator function as a mathematical tool for examining the geometric structures of designs. In Section 3 geometric isomorphism is defined for the classification of factorial designs. Section 4 proposes a new aberration criterion for factorial designs with quantitative factors. Some remarks are given in Section 5.

In the remainder of this section we will introduce some notation and terminology. Let $\mathcal{D}$ be the $\mathrm{OA}(N, s_1 s_2 \ldots s_k)$, which is a full factorial design with $k$ factors and $N$ design points, where $N = s_1 s_2 \ldots s_k$. Unless specified, the levels of $i$th factor are set at $G_i = \{0, 1, \ldots, s_i - 1\} \subset \mathbb{R}$ for each factor, which are evenly spaced. Therefore, $\mathcal{D}$ is a set of $N$ points in $\mathbb{R}^k$. A $k$-factor factorial design $\mathcal{A}$ is said to be in a *design space* $\mathcal{D}$ if its design points are all in $\mathcal{D}$, that is, $\forall \mathbf{x} \in \mathcal{A}, \mathbf{x} \in \mathcal{D}$. A design point in $\mathcal{D}$ may appear more than once in $\mathcal{A}$. Throughout this paper, $\sum_{\mathbf{x} \in \mathcal{A}} f(\mathbf{x})$ sums the function $f$ over all design points in $\mathcal{A}$; that is, if $\mathbf{x}$ appears multiple times, $f(\mathbf{x})$ is summed over multiple times.

For each factor $X_i$, define a set of orthogonal contrasts $C_0^i(x), C_1^i(x), \ldots, C_{s_i-1}^i(x)$ such that

$$\sum_{x \in \{0,1,\ldots,s_i-1\}} C_u^i(x) C_v^i(x) = \begin{cases} 0, & \text{if } u \neq v, \\ s_i, & \text{if } u = v. \end{cases} \tag{1.1}$$

Let $\mathcal{T} = G_1 \times \cdots \times G_k$. An *orthonormal contrast basis* (OCB) on $\mathcal{D}$ is defined as

$$C_\mathbf{t}(\mathbf{x}) = \prod_{i=1}^{k} C_{t_i}^i(x_i) \tag{1.2}$$

for $\mathbf{t} = (t_1, t_2, \ldots, t_k) \in \mathcal{T}$ and $\mathbf{x} = (x_1, x_2, \ldots, x_k) \in \mathcal{D}$. It is obvious that

$$\sum_{\mathbf{x} \in \mathcal{D}} C_\mathbf{t}(\mathbf{x}) C_\mathbf{u}(\mathbf{x}) = \begin{cases} 0, & \text{if } \mathbf{t} \neq \mathbf{u}, \\ N, & \text{if } \mathbf{t} = \mathbf{u}, \end{cases} \tag{1.3}$$

where $\mathbf{t}, \mathbf{u}$ are elements in the set $\mathcal{T}$. In statistical analysis, $C_0^i(x) = 1$ is often adopted to represent a constant term. Therefore, we call $\{C_\mathbf{t}(\mathbf{x})\}$ with $C_0^i = 1$ for all $i$ a *statistical orthonormal contrast basis* (SOCB). When $C_j^i(x)$ is a polynomial of degree $j$ for $j = 0, 1, \ldots, s_i - 1$ and $i = 1, 2, \ldots, k$, the SOCB is called an *orthogonal polynomial basis* (OPB) [Draper and Smith (1998), Chapter 22]. Note that an OPB is an SOCB, and an SOCB is an OCB.

We define two norms on $\mathcal{T}$. Let $\|\mathbf{t}\|_0$ be the number of nonzero elements in $\mathbf{t}$ and let

$$\|\mathbf{t}\|_1 = \sum_{i=1}^{k} t_i.$$

For a contrast $C_\mathbf{t}$ in an SOCB, $\|\mathbf{t}\|_0$ is the number of factors it involves. If the SOCB is also an OPB, $\|\mathbf{t}\|_1$ gives its polynomial degree. Two contrasts $C_\mathbf{t}$ and $C_\mathbf{u}$ in an SOCB are said to be *disjoint* if they have no common factors, that is, $\max_{1 \leq i \leq k} \min(t_i, u_i) = 0$.

**2. Indicator functions.** Indicator functions are presented in Fontana, Pistone and Rogantin (2000) for studying two-level fractional factorial designs (without replicates). Ye (2003) generalizes to accommodate replicates. In this section the definition is extended further to general factorial designs.

DEFINITION 2.1. Let $\mathcal{A}$ be a design in the design space $\mathcal{D}$. The indicator function $F_\mathcal{A}(\mathbf{x})$ of $\mathcal{A}$ is a function defined on $\mathcal{D}$, such that for $\mathbf{x} \in \mathcal{D}$, the value of $F_\mathcal{A}(\mathbf{x})$ is the number of appearances of point $\mathbf{x}$ in design $\mathcal{A}$.

The following proposition follows immediately from the definition.

PROPOSITION 2.1. *Let $\mathcal{A}_1, \ldots, \mathcal{A}_m$ be factorial designs of the same design space $\mathcal{D}$ and $F_{\mathcal{A}_i}(\mathbf{x})$, $i = 1, \ldots, m$, be their corresponding indicator functions. Let $\mathcal{B}$ be the combined design (design points are repeatable in $\mathcal{B}$) of $\mathcal{A}_1, \ldots, \mathcal{A}_m$. Then the indicator function of $\mathcal{B}$ is*

$$F_{\mathcal{B}}(\mathbf{x}) = \sum_{i=1}^{m} F_{\mathcal{A}_i}(\mathbf{x}).$$

Since a design is uniquely represented by its indicator function, the indicator function carries all properties of this design. Some of these properties are revealed when indicator functions are expanded with respect to an OCB.

THEOREM 2.1. *Let $\mathcal{A}$ be a factorial design with $n$ runs. Let $\mathcal{D}$ be the design space of $\mathcal{A}$, and $\{C_{\mathbf{t}}(\mathbf{x}), \mathbf{t} \in \mathcal{T}\}$ be an OCB defined on $\mathcal{D}$. The indicator function of $\mathcal{A}$ can be represented as a linear combination of $C_{\mathbf{t}}s$ as follows:*

$$(2.1) \qquad F_{\mathcal{A}}(\mathbf{x}) = \sum_{\mathbf{t} \in \mathcal{T}} b_{\mathbf{t}} C_{\mathbf{t}}(\mathbf{x}),$$

*for all $\mathbf{x} \in \mathcal{D}$. The coefficients $\{b_{\mathbf{t}}, \mathbf{t} \in \mathcal{T}\}$ are uniquely determined as*

$$(2.2) \qquad b_{\mathbf{t}} = \frac{1}{N} \sum_{\mathbf{x} \in \mathcal{A}} C_{\mathbf{t}}(\mathbf{x}).$$

*And, in particular, for an SOCB, $b_{\mathbf{0}} = n/N$, where $\mathbf{0} = (0, 0, \ldots, 0)$.*

PROOF. The indicator function $F_{\mathcal{A}}(\mathbf{x})$ is defined on $\mathcal{D}$ and $F_{\mathcal{A}}(\mathcal{D})$ can be viewed as a vector in $\mathbb{R}^N$. Since the $\{C_{\mathbf{t}}(\mathcal{D}), \mathbf{t} \in \mathcal{T}\}$ forms a basis of $\mathbb{R}^N$, $F_{\mathcal{A}}(\mathcal{D})$ can be represented as a linear combination of $\{C_{\mathbf{t}}(\mathcal{D})\}$. Equivalently, (2.1) is true. For the coefficients $b_{\mathbf{t}}s$,

$$\sum_{\mathbf{x} \in \mathcal{A}} C_{\mathbf{t}}(\mathbf{x}) = \sum_{\mathbf{x} \in \mathcal{D}} F_{\mathcal{A}}(\mathbf{x}) C_{\mathbf{t}}(\mathbf{x}) = \sum_{\mathbf{x} \in \mathcal{D}} \sum_{\mathbf{s} \in \mathcal{T}} b_{\mathbf{s}} C_{\mathbf{s}}(\mathbf{x}) C_{\mathbf{t}}(\mathbf{x})$$
$$= \sum_{\mathbf{s} \in \mathcal{T}} b_{\mathbf{s}} \sum_{\mathbf{x} \in \mathcal{D}} C_{\mathbf{s}}(\mathbf{x}) C_{\mathbf{t}}(\mathbf{x}) = N b_{\mathbf{t}}.$$

The proof is complete. □

Note that the theorem does not depend on level settings and choice of $C_j^i(x)$, as long as (1.1) is satisfied and $\{C_{\mathbf{t}}\}$ is defined as in (1.2). In the functional space generated by linear combinations of $\{C_{\mathbf{t}}\}$, the indicator function of a design has a unique representation, that is, there is a one-to-one relation between a factorial design and its $b_{\mathbf{t}}$ values. This is an extension of a similar result on two-level designs presented in Ye (2003). When $\{C_{\mathbf{t}}\}$ is an OPB, an indicator function can be uniquely represented as a polynomial of degree no more than $\prod_{i=1}^{k}(s_i - 1)$.

A projected design has the same number of runs as the original design but is in a reduced design space with only a subset of the original factors. Given a design's polynomial representation in the form of (2.1), the polynomial representations of its projected designs are easily available, as shown in the following corollary.

COROLLARY 2.1. *Let $\mathcal{A}$ be a factorial design in design space $\mathcal{D}$ and $F_{\mathcal{A}}(\mathbf{x}) = \sum_{t \in \mathcal{T}} b_{\mathbf{t}} C_{\mathbf{t}}(\mathbf{x})$ be its indicator function. Without loss of generality, let $\mathcal{B}$ be its projection to factors $X_1, \ldots, X_l$. If $\{C_{\mathbf{t}}\}$ is an SOCB, the indicator function of $\mathcal{B}$ is then*

$$F_{\mathcal{B}}(x_1, \ldots, x_l) = N_2 \sum_{\mathbf{t} \in \mathcal{T}_1} b_{\mathbf{t}} C_{\mathbf{t}}, \tag{2.3}$$

*where*

$$N_2 = \prod_{i=l+1}^{k} s_i \quad \text{and} \quad \mathcal{T}_1 = \{\mathbf{t} | t_{l+1} = \cdots = t_k = 0\}.$$

PROOF. From (2.2), $b_{\mathbf{t}} = 1/N \sum_{\mathbf{x} \in \mathcal{A}} C_{\mathbf{t}}(\mathbf{x})$. The coefficient of $C_{\mathbf{t}}(\mathbf{x})$ in $F_{\mathcal{B}}(\mathbf{x})$ is then $1/N_1 \sum_{\mathbf{x} \in \mathcal{B}} C_{\mathbf{t}}(\mathbf{x})$, where $N_1 = \prod_{i=1}^{l} s_i$. Equation (2.3) follows. □

The coefficients $b_{\mathbf{t}}$ also relate to the orthogonality of a design. This can be shown in the following corollary which follows immediately from (1.2) and (2.2).

COROLLARY 2.2. *Let $\{C_{\mathbf{t}}(\mathbf{x}), \mathbf{t} \in \mathcal{T}\}$ be an SOCB. For disjoint $C_{\mathbf{u}}$ and $C_{\mathbf{v}}$,*

$$b_{\mathbf{u}+\mathbf{v}} = \frac{1}{N} \sum_{\mathbf{x} \in \mathcal{A}} C_{\mathbf{u}}(\mathbf{x}) C_{\mathbf{v}}(\mathbf{x}).$$

*Furthermore, the correlation of $C_{\mathbf{u}}$ and $C_{\mathbf{v}}$ in $\mathcal{A}$ is $b_{\mathbf{u}+\mathbf{v}}/b_{\mathbf{0}}$.*

Let $C_{\mathbf{u}}(\mathbf{x})$ and $C_{\mathbf{v}}(\mathbf{x})$ be two disjoint contrasts. From Corollary 2.2, the two contrasts are zero correlated on design $\mathcal{A}$ if and only if $b_{\mathbf{u}+\mathbf{v}} = 0$. As a special case, $b_{\mathbf{t}} = 0$ implies that the contrast $C_{\mathbf{t}}(\mathbf{x})$ has zero correlation with the constant term on design $\mathcal{A}$. In general, a smaller $b_{\mathbf{t}}$ implies a lesser degree of aliasing between effects and, therefore, $b_{\mathbf{t}}$ can be used as a measurement of aliasing between effects. Various statistical properties of designs can be studied through the $b_{\mathbf{t}}$'s. More results will be shown in Sections 3 and 4.

EXAMPLE 2.1. Consider the case of three-level factorial designs with $k$ factors. The design space $\mathcal{D}$ is a collection of $3^k$ points: $\{(d_1, \ldots, d_k), d_i = 0, 1, 2, i = 1, \ldots, k\}$. From Definition 2.1, any $k$-factor three-level factorial design can be represented by an indicator function defined on $\mathcal{D}$. The orthonormal polynomials for a three-level factor are

$$C_0(x) = 1, \qquad C_1(x) = \sqrt{\tfrac{3}{2}}(x-1) \quad \text{and} \quad C_2(x) = \sqrt{2}(\tfrac{3}{2}(x-1)^2 - 1).$$

Note that $(C_1(0), C_1(1), C_1(2)) = (-\sqrt{3/2}, 0, \sqrt{3/2})$ and $(C_2(0), C_2(1), C_2(2)) = (1/\sqrt{2}, -\sqrt{2}, 1/\sqrt{2})$ are proportional to the linear and quadratic contrasts, respectively, as defined in Wu and Hamada [(2000), Section 5.6]. Thus, $\{C_{\mathbf{t}}(\mathbf{x}), \mathbf{t} \in \mathcal{T}\}$, where $\mathcal{T}$ is the vector space $\{0, 1, 2\}^k$, is an OPB for the functional space of $\mathcal{D}$. By Theorem 2.1, an indicator function can be written as a linear combination of $C_{\mathbf{t}}(\mathbf{x})$'s with coefficients

$$b_{\mathbf{t}} = \frac{1}{3^k} \sum_{\mathbf{x} \in \mathcal{A}} C_{\mathbf{t}}(\mathbf{x})$$

and, in particular, $b_{\mathbf{0}} = n/3^k$. The coefficients $b_{\mathbf{t}}$ contain information about aliasing between effects. For the design on the right-hand side of Table 1, which is also shown on the right-hand side of Figure 1, its indicator function is

$$\begin{aligned}(2.4) \quad F(\mathbf{x}) = {}& \tfrac{1}{3}C_{000}(\mathbf{x}) - \tfrac{\sqrt{6}}{12}C_{111}(\mathbf{x}) - \tfrac{\sqrt{2}}{12}C_{112}(\mathbf{x}) + \tfrac{\sqrt{2}}{12}C_{121}(\mathbf{x}) + \tfrac{\sqrt{2}}{12}C_{211}(\mathbf{x}) \\ & - \tfrac{\sqrt{6}}{12}C_{122}(\mathbf{x}) - \tfrac{\sqrt{6}}{12}C_{212}(\mathbf{x}) + \tfrac{\sqrt{6}}{12}C_{221}(\mathbf{x}) + \tfrac{\sqrt{2}}{12}C_{222}(\mathbf{x}),\end{aligned}$$

where

$$C_{i_1 i_2 i_3}(\mathbf{x}) = C_{i_1}(x_1) C_{i_2}(x_2) C_{i_3}(x_3).$$

For the design on the left-hand side, its indicator function is

$$(2.5) \quad F(\mathbf{x}) = \tfrac{1}{3}C_{000}(\mathbf{x}) + \tfrac{\sqrt{2}}{6}C_{112}(\mathbf{x}) + \tfrac{\sqrt{2}}{6}C_{121}(\mathbf{x}) + \tfrac{\sqrt{2}}{6}C_{211}(\mathbf{x}) - \tfrac{\sqrt{2}}{6}C_{222}(\mathbf{x}).$$

In (2.5), $b_{111} = 0$ implies that the linear-by-linear interactions are orthogonal to linear main effects in the design. For the other design, they are not orthogonal since $b_{111} = \tfrac{\sqrt{6}}{12}$ in (2.4). This is consistent with the higher $D$-efficiency of the former design when a model with all linear main effects and linear-by-linear interactions is considered.

**3. Geometric isomorphism.** When factors are all nominal in a factorial design, new design matrices obtained through level permutations in one or more factors are considered to be isomorphic to the original design. This is referred to as the combinatorial isomorphism. As shown in Section 1, when factors are quantitative, level permutations generate designs with different geometric structures and, thus, different design properties. Cheng and Wu (2001) observed differences in $D$-efficiency of these designs and proposed model isomorphism for classification. However, such classification depends on *a priori* specified models. Designs that have the same efficiencies with respect to a given model might have different efficiencies with respect to another model. A classification with respect to a certain model can be no longer useful when a different model is considered. For consistency, a better classification method should be based on the geometric structures, which are fundamental to design properties and do not depend on the choice of the models.

From a geometric viewpoint, a geometric object remains the same structure when rotating and/or reflecting with respect to a super-plane. In the context of a factorial design, only rotations and reflections, after which the resulting designs are still in the design space $\mathcal{D}$, should be considered. Rotating then corresponds to variable exchange and reflecting corresponds to reversing order of the levels. Therefore, we define *geometric isomorphism* of two designs as follows.

DEFINITION 3.1. Let $\mathcal{A}$ and $\mathcal{B}$ be two factorial designs from the same design space $\mathcal{D}$. Designs $\mathcal{A}$ and $\mathcal{B}$ are said to be geometrically isomorphic if one can be obtained from the other by variable exchange and/or reversing the level order of one or more factors.

One can differentiate geometrically nonisomorphic designs by comparing their indicator functions. Let $F_{\mathcal{A}}(x_1, \ldots, x_k)$ be the indicator function of design $\mathcal{A}$ and $\mathcal{A}_1$ be the design obtained by reversing the level order of factor $X_1$ in $\mathcal{A}$. The indicator function of $\mathcal{A}_1$ is $F_{\mathcal{A}}(2d - x_1, x_2, \ldots, x_k)$, where $d$ is the center of the levels. Let $\mathcal{A}_2$ be the design obtained by exchanging factor $X_1$ with factor $X_2$ in $\mathcal{A}$. Then its indicator function is $F_{\mathcal{A}}(x_2, x_1, \ldots, x_k)$. If the indicator functions of two designs are the same after a series of such operations, then they are geometrically isomorphic.

Geometric isomorphism of two designs can be more easily examined when indicator functions are expanded to polynomial form with respect to an OPB. Theorem 3.1 implies that if two designs are geometrically isomorphic, the absolute values of their coefficients $b_\mathbf{t}$ must show the same frequency patterns.

THEOREM 3.1. *Let $\mathcal{A}$ and $\mathcal{B}$ be two factorial designs of the design space $\mathcal{D}$, and $\{C_\mathbf{t}(\mathbf{x})\}$ be an OPB defined on $\mathcal{D}$. Let $F_{\mathcal{A}}(\mathbf{x}) = \sum b_\mathbf{t} C_\mathbf{t}(\mathbf{x})$ and $F_{\mathcal{B}}(\mathbf{x}) = \sum b'_\mathbf{t} C_\mathbf{t}(\mathbf{x})$ be the indicator functions of $\mathcal{A}$ and $\mathcal{B}$, respectively. Designs $\mathcal{A}$ and $\mathcal{B}$ are geometrically isomorphic if and only if there exist a* permutation $(i_1 i_2 \ldots i_k)$ and a vector $(j_1 j_2 \ldots j_k)$, where $j_l$'s are either 0 or 1, such that

$$(3.1) \qquad b_{t_1 t_2 \ldots t_k} = \left( \prod_{l=1}^{k} (-1)^{j_l t_l} \right) b'_{t_{i_1} t_{i_2} \ldots t_{i_k}}$$

for all $\mathbf{t} = (t_1 t_2 \ldots t_k) \in \mathcal{T}$.

PROOF. Using the three-term recursive equation given in Kennedy and Gentle [(1980), pages 343 and 344] for constructing orthogonal polynomials, it is easy to show that the orthogonal polynomial contrasts $C_j(x)$ of a factor satisfy the following condition:

$$(3.2) \qquad C_j(x) = \begin{cases} -C_j(2d - x), & \text{if } j \text{ odd}, \\ C_j(2d - x), & \text{if } j \text{ even}. \end{cases}$$

If $\mathcal{A}$ and $\mathcal{B}$ are geometrically isomorphic, then by definition $\mathcal{A}$ must be obtained from $\mathcal{B}$ by variable exchange and/or reversal of levels. Let the variable exchange be $x_l \to x_{i_l}$, and let $j_l = 1$ if the levels of factor $x_l$ are reversed, $j_l = 0$ if not. Hence, (3.1) is truly based on (2.2) and (3.2). Conversely, if (3.1) is true, $\mathcal{B}$ can be obtained from $\mathcal{A}$ by the variable exchange $x_l \to x_{i_l}$ and the level reverses on the factors with $j_l = 1$. Therefore, $\mathcal{A}$ and $\mathcal{B}$ are geometrically isomorphic. □

Note that from the proof Theorem 3.1, it holds for any basis such that (3.2) is satisfied. From the theorem, one can immediately show that (2.4) and (2.5) represent two geometrically nonisomorphic designs as their coefficients show different frequency patterns. In general, with a proper choice of $\{C_\mathbf{t}(\mathbf{x})\}$, two designs are geometrically isomorphic if and only if their indicator functions have the same coefficients $b_\mathbf{t}$ after a certain type of permutation and sign reversal. Otherwise, if two designs have different geometric structures, their coefficients must show different frequency patterns.

EXAMPLE 3.1. The $L_{18}$ array (shown in Table 2) is one of the most popular designs among industrial experimenters. It has one two-level factor and seven three-level factors. For the moment, we only consider the three-level factors. Wang and Wu (1995) studied the projected design of $L_{18}$ and reported three combinatorially nonisomorphic cases for 3-factor projections (denoted as 18-3.1, 18-3.2 and 18-3.3) and four combinatorially nonisomorphic cases for 4-factor projections (denoted as 18-4.1, 18-4.2, 18-4.3, 18-4.4). As shown earlier in this paper, level permutation in a design may create designs with different geometric structures. There are a total of six permutations among three levels, that is,

$\{0, 1, 2\} \to \{0, 1, 2\}, \qquad \{0, 1, 2\} \to \{0, 2, 1\}, \qquad \{0, 1, 2\} \to \{1, 0, 2\},$

$\{0, 1, 2\} \to \{1, 2, 0\}, \qquad \{0, 1, 2\} \to \{2, 0, 1\} \quad \text{and} \quad \{0, 1, 2\} \to \{2, 1, 0\}.$

The six permutations can be divided into three pairs as shown in Table 3. Within each pair, one permutation is the reverse of the other, hence, only one is needed in generating geometrically nonisomorphic designs. For each combinatorially nonisomorphic case, permutations are applied to each column to search for all geometrically nonisomorphic cases. For three-factor projections, there are two, four and two geometrically nonisomorphic cases within 18-3.1, 18-3.2 and 18-3.3, respectively. Cheng and Wu (2001) reported the same number of model nonisomorphic cases for 18-3.1 and 18-3.2 but did not report model nonisomorphic cases for 18-3.3. For four-factor projection, there are four, ten, three and four geometrically nonisomorphic cases in 18.4-1, 18.4-2, 18-4.3 and 18-4.4, respectively. Cheng and Wu (2001) only reported four model nonisomorphic cases for 18-4.2 and none for the other three. A complete list of these geometrically nonisomorphic projected designs is given in the Appendix.

**4. Aberration criterion.** A popular criteria for factorial designs is minimum aberration. The original definition of minimum aberration based on group theory applies to regular $p^{m-n}$ fractional factorial designs [Fries and Hunter (1980)]. Recently, Xu and Wu (2001) proposed an aberration criterion based on coding theory for general factorial designs. It reduces to the traditional aberration criterion for regular $p^{m-n}$ designs, and the $G_2$ aberration criteria [Tang and Deng (1999)] for general two-level factorial designs. A statistical justification is given by Xu and Wu (2001) to relate the criterion with ANOVA. From this relation, it can be easily seen that their aberration criterion can be redefined using the indicator functions as follows.

Table 2
$L_{18}$ *orthogonal array*

| 0 | 1 | 2 | 3 | 4 | 5 | 6 | 7 |
|---|---|---|---|---|---|---|---|
| 0 | 0 | 0 | 0 | 0 | 0 | 0 | 0 |
| 0 | 0 | 1 | 1 | 1 | 1 | 1 | 1 |
| 0 | 0 | 2 | 2 | 2 | 2 | 2 | 2 |
| 0 | 1 | 0 | 0 | 1 | 1 | 2 | 2 |
| 0 | 1 | 1 | 1 | 2 | 2 | 0 | 0 |
| 0 | 1 | 2 | 2 | 0 | 0 | 1 | 1 |
| 0 | 2 | 0 | 1 | 0 | 2 | 1 | 2 |
| 0 | 2 | 1 | 2 | 1 | 0 | 2 | 0 |
| 0 | 2 | 2 | 0 | 2 | 1 | 0 | 1 |
| 1 | 0 | 0 | 2 | 2 | 1 | 1 | 0 |
| 1 | 0 | 1 | 0 | 0 | 2 | 2 | 1 |
| 1 | 0 | 2 | 1 | 1 | 0 | 0 | 2 |
| 1 | 1 | 0 | 1 | 2 | 0 | 2 | 1 |
| 1 | 1 | 1 | 2 | 0 | 1 | 0 | 2 |
| 1 | 1 | 2 | 0 | 1 | 2 | 1 | 0 |
| 1 | 2 | 0 | 2 | 1 | 2 | 0 | 1 |
| 1 | 2 | 1 | 0 | 2 | 0 | 1 | 2 |
| 1 | 2 | 2 | 1 | 0 | 1 | 2 | 0 |

Table 3
*Six permutations of three levels*

| $I$ | $u$ | $u^2$ |
|---|---|---|
| 0 | 1 | 2 |
| 1 | 2 | 0 |
| 2 | 0 | 1 |
| ↻ | ↻ | ↻ |
| 2 | 1 | 0 |
| 1 | 0 | 2 |
| 0 | 2 | 1 |

DEFINITION 4.1. Let $\mathcal{A}$ be an $n \times k$ factorial design of design space $\mathcal{D}$. Let $F_\mathcal{A}(\mathbf{x}) = \sum_{t \in \mathcal{T}} b_\mathbf{t} C_\mathbf{t}(\mathbf{x})$ be its indicator function, where $\{C_\mathbf{t}\}$ is an SOCB. The generalized wordlength pattern $(\alpha_1(\mathcal{A}), \ldots, \alpha_k(\mathcal{A}))$ of design $\mathcal{A}$ is defined as

$$(4.1) \qquad \alpha_i(\mathcal{A}) = \sum_{\|\mathbf{t}\|_0 = i} \left(\frac{b_\mathbf{t}}{b_\mathbf{0}}\right)^2.$$

The generalized minimum aberration criterion is to sequentially minimize $\alpha_i(\mathcal{A})$ for $i = 1, 2, \ldots, k$. The resolution of $\mathcal{A}$ equals the smallest $r$ such that $\alpha_r > 0$.

From Corollary 2.2, $(b_\mathbf{t}/b_\mathbf{0})^2$ is a measurement that reflects the severity of aliasing between the effect $C_\mathbf{t}$ and the general mean. Therefore, in (4.1), $\alpha_i$ measures the overall aliasing between all $i$-factor effects and the general mean. A smaller $\alpha_i$ indicates a lesser degree of aliasing between the $i$-factor effects and the overall mean. Therefore, the $\alpha_i$'s are to be minimized sequentially. Note that the definition implicitly assumes all $i$-factor effects are equally important, which is only suitable for nominal factors (see later discussion on the hierarchical ordering principle). The above definition is a natural generalization of the definition of the aberration criterion for two-level factorial designs given in Ye (2003). Xu and Wu (2001) showed that $\mathcal{A}$ is an orthogonal array of strength $t$ if and only if $\alpha_i(\mathcal{A}) = 0$ for $1 \leq i \leq t$. From the definition, the if and only if condition can be stated in the language of indicator function as follows: $b_{\mathbf{t}} = 0$ for all $\mathbf{t}$'s such that $1 \leq \|\mathbf{t}\|_0 \leq t$. For example, the coefficients $b_{\mathbf{t}}$ in indicator functions (2.4) and (2.5) are zero for all $\mathbf{t}$ such that $1 \leq \|\mathbf{t}\|_0 \leq 2$. Therefore, they are orthogonal arrays of strength two. This aberration criterion is invariant to level permutation as well as the choice of contrasts. While these features are desirable when all factors are nominal, it is not quite desirable when factors are quantitative. An immediate problem is that this aberration criterion completely fails to distinguish and rank combinatorially isomorphic but geometrically nonisomorphic designs, for example, the two designs in Table 1.

An important assumption behind aberration criteria is the *hierarchical ordering principle* [Wu and Hamada (2000), Section 3.5]: (i) low-order effects are more likely to be important than high-order effects, and (ii) effects of the same order are equally likely to be important. The principle can be applied to nominal and quantitative factors. However, the effect orders for the two types of factors should be different. For nominal factors, the objective of analysis is treatment comparison. Therefore, all $i$-factor effects are regarded as equally important and $i$-factor effects are more important than $j$-factor effects for $i < j$. The effect order is decided by the number of factors that are related to the corresponding contrast $C_{\mathbf{t}}$, that is, the value of $\|\mathbf{t}\|_0$. Therefore, in (4.1), the overall aliasing is measured by taking the sum over those $\mathbf{t}$'s with the same $\|\mathbf{t}\|_0$ value. For experiments with quantitative factors, polynomial models are often utilized to approximate the response. In this case, effects of higher polynomial degree are regarded as less important than effects of lower polynomial degree. Therefore, the order of effect importance should be arranged according to polynomial degrees. Recall that in an OPB, $C_j(x)$ is a polynomial of degree $j$ and $C_{\mathbf{t}}(x)$ is a polynomial of degree $\|\mathbf{t}\|_1$. In this case, the order of effect importance can be defined according to the values of $\|\mathbf{t}\|_1$. For example, for quantitative three-level factors, the order of effect importance is as follows:

$$
\begin{aligned}
&l >> q == ll >> lq == ql == lll \\
&>> qq == llq == lql == qll == llll >> \cdots,
\end{aligned}
\tag{4.2}
$$

where $>>$ is read as "more important than" and $==$ as "as important as," and $l$ and $q$ indicate linear and quadratic main effects, respectively, $ll$ linear-by-linear interaction, etc. Such ordering is consistent with the response surface methodology in which, based on a rationale from Taylor's series expansion, effects of the same degree are sequentially added to the model starting from the lowest degree.

Consider the two designs in Table 1. When all factors are quantitative and OPB is used, the contrasts in equations (2.4) and (2.5) follow the linear-quadratic decomposition and have clear interpretation in terms of fitting polynomial models. As mentioned previously, in (2.5), $b_{111} = 0$ implies that the *lll* interactions have zero correlation with the constant term, and the *ll* interactions between any two factors have zero correlation with the linear main effect of the third factor. However, the design of (2.4) does not have this nice property since $b_{111} = \frac{\sqrt{6}}{12}$; hence, the *lll* interaction is (partially) aliased with the constant term, similar with the linear main effects and *ll* interactions. In Xu and Wu (2001), $b_{111}$ is considered as important as $b_{112}, \ldots, b_{222}$ and are summed up together in $\alpha_3$. This is not quite appropriate if a polynomial model is to be fitted and *lll* interactions are regarded as more important than all other three-way interactions as in (4.2). Therefore, when all factors are quantitative with polynomial models as our point of interest, the following aberration and resolution criteria are proposed.

DEFINITION 4.2. Let $\mathcal{A}$ be an $n \times k$ factorial design with quantitative factors, and let $\{C_\mathbf{t}\}$ be an OPB. Let $F_\mathcal{A}(\mathbf{x}) = \sum_{t \in \mathcal{T}} b_\mathbf{t} C_\mathbf{t}(\mathbf{x})$ be the indicator function of $\mathcal{A}$. The generalized wordlength pattern $(\beta_1(\mathcal{A}), \ldots, \beta_K(\mathcal{A}))$ is defined as

$$\beta_i(\mathcal{A}) = \sum_{\|\mathbf{t}\|_1 = i} \left(\frac{b_\mathbf{t}}{b_\mathbf{0}}\right)^2. \tag{4.3}$$

The generalized minimum aberration criterion is to sequentially minimize $\beta_i$ for $i = 1, 2, \ldots, K$, where

$$K = \sum_{i=1}^k (s_i - 1)$$

is the highest possible degree in the decomposition. The resolution of $\mathcal{A}$ is defined to be the smallest $r$ such that $\beta_r > 0$.

In the rest of this paper we refer to the word length pattern given in Definition 4.1 as the $\alpha$ wordlength pattern, and the wordlength pattern in the above definition as the $\beta$ wordlength pattern. It is probably more appropriate to call the new aberration and resolution criteria "*polynomial degree*" aberration and resolution. For three-level designs, (4.3) counts the overall aliasing between the general mean and effects that are of the same importance in (4.2). Based on the above definition, the wordlength patterns of the two designs in (2.4) and (2.5) are $(0, 0, \frac{3}{8}, \frac{3}{8}, \frac{9}{8}, \frac{1}{8})$ and $(0, 0, 0, \frac{3}{2}, 0, \frac{1}{2})$, respectively. The latter one has less aberration and higher resolution and is favored.

Two contrasts, $C_\mathbf{t}(x)$ and $C_\mathbf{u}(x) \in \mathcal{T}$, are treated as equally important if and only if $\|\mathbf{t}\|_1 = \|\mathbf{u}\|_1$. Therefore, by Theorem 3.1, we can easily show that geometrically isomorphic designs have identical $\beta$ wordlength patterns.

COROLLARY 4.1. *Let $\mathcal{A}$ and $\mathcal{B}$ be two geometrically isomorphic designs in design space $\mathcal{D}$. Then their $\beta$ wordlength patterns are identical.*

PROOF. Let the indicator functions of two designs be $F_{\mathcal{A}} = \sum b_{\mathbf{t}} C_{\mathbf{t}}(x)$ and $F_{\mathcal{B}} = \sum b'_{\mathbf{t}} C_{\mathbf{t}}(x)$. By Theorem 3.1, there must be a permutation $(i_1 i_2 \ldots i_k)$ such that $(b_{t_1 t_2 \ldots t_k})^2 = (b'_{t_{i_1} t_{i_2} \ldots t_{i_k}})^2$ for all $\mathbf{t} = (t_1 t_2 \ldots t_k) \in \mathcal{T}$. For each $\mathbf{t}$, denote its corresponding permutation as $\mathbf{t}'$. Since $\|\mathbf{t}\|_1 = \|\mathbf{t}'\|_1$, $\beta_j(\mathcal{A})$ and $\beta_j(\mathcal{B})$ sum over the same values, and hence, are identical. $\square$

By Corollary 4.1, two designs with different $\beta$ wordlength patterns must be geometrically nonisomorphic.

The only distinction between Definitions 4.1 and 4.2 are the norms of $\mathbf{t}$ used in (4.1) and (4.3), which reflect the difference in ordering effects for nominal and quantitative factors. For a given design, the sum of its $\alpha_i$'s is the same as the sum of its $\beta_i$'s. The following theorem shows that this sum is a constant for designs that have the same run sizes and replication patterns.

THEOREM 4.1. *Let $\mathcal{A}$ be an $n \times k$ factorial design in the design space $\mathcal{D}$. Let $F_{\mathcal{A}}(\mathbf{x}) = \sum_{\mathbf{t} \in \mathcal{T}} b_{\mathbf{t}} C_{\mathbf{t}}(\mathbf{x})$ be the indicator function of $\mathcal{A}$. Then*

$$\text{(4.4)} \qquad \sum_{\mathbf{t} \in \mathcal{T}} \left( \frac{b_{\mathbf{t}}}{b_{\mathbf{0}}} \right)^2 = \frac{n_2 N}{n^2},$$

*where*

$$n_2 = \sum_{\mathbf{x} \in \mathcal{D}} F_{\mathcal{A}}^2(\mathbf{x}) \quad \text{and} \quad N = s_1 \ldots s_k.$$

*For designs with no replicates, $n_2 = n$.*

PROOF.

$$\sum_{\mathbf{x} \in \mathcal{D}} F_{\mathcal{A}}^2(\mathbf{x}) = \sum_{\mathbf{x} \in \mathcal{D}} \left( \sum_{\mathbf{t} \in \mathcal{T}} b_{\mathbf{t}} C_{\mathbf{t}} \right)^2 = \sum_{\mathbf{x} \in \mathcal{D}} \sum_{\mathbf{t}_1, \mathbf{t}_2 \in \mathcal{T}} b_{\mathbf{t}_1} b_{\mathbf{t}_2} C_{\mathbf{t}_1}(\mathbf{x}) C_{\mathbf{t}_2}(\mathbf{x})$$

$$= \sum_{\mathbf{t}_1, \mathbf{t}_2 \in \mathcal{T}} b_{\mathbf{t}_1} b_{\mathbf{t}_2} \sum_{\mathbf{x} \in \mathcal{D}} C_{\mathbf{t}_1}(\mathbf{x}) C_{\mathbf{t}_2}(\mathbf{x}) = N \sum_{\mathbf{t} \in \mathcal{T}} b_{\mathbf{t}}^2.$$

From Theorem 2.1, $b_{\mathbf{0}} = n/N$. Hence, (4.4) is obtained. For designs with no replicates, $F_{\mathcal{A}}^2(\mathbf{x}) = F_{\mathcal{A}}(\mathbf{x})$, hence $n_2 = n$. $\square$

A special case of the above theorem is a well-known result for regular fractional factorial $p^{k-m}$ designs, in which the sum of their wordlength pattern vector equals $p^m - 1$. The theorem shows that it holds for both $\alpha$ and $\beta$ wordlength patterns. The theorem also shows that the sum of wordlength pattern vectors is larger for designs with higher degrees of replication as $n_2$ is larger in (4.4). Therefore, they tend to have higher aberration than those with less replicates.

TABLE 4
*Minimum aberration projections of $L_{18}$ with only three-level factors*

| # of factors | Columns | $(\beta_3, \beta_4, \beta_5)$ | Resolution |
|---|---|---|---|
| 3 | $1u^2$ 2 5 | (0, 0.125, 0.75) | IV |
| 4 | $1u^2$ $2u^2$ $3u^2$ 5 | (0, 1.875, 0) | IV |
| 5 | $1u^2$ $2u^2$ $3u^2$ 4 5 | (0, 6.0625, 0) | IV |
| 6 | $2u$ $3u$ $4u$ 5 6 7 | (0.75, 6.9375, 6.75) | III |
| 7 | $1u$ $2u$ $3u$ $4u$ 5 6 7 | (1.5, 14.625, 12) | III |

Note that for two-level designs, both $\alpha$ and $\beta$ wordlength patterns reduce to the same generalized wordlength pattern by Tang and Deng (1999) and Ye (2003). However, if the factors have more than two levels, these two wordlength patterns often give different minimum aberration designs. One should choose from them based on the nature of the factors, nominal or quantitative.

EXAMPLE 4.1. When less than seven three-level factors are considered in an experiment, it would be of interest to know which columns in the $L_{18}$ array are the best to be assigned to those factors. To find the minimum aberration projections of the $L_{18}$ array, an exhaustive search over all possible projections was performed, and three-level permutations were applied to each column. Tables 4 and 5 list the $\beta$ minimum aberration projected designs with and without the two-level factors, respectively. In the tables, $u$ denotes the permutation $\{0, 1, 2\} \to \{1, 2, 0\}$; $u^2$ denotes the permutation $\{0, 1, 2\} \to \{2, 0, 1\}$. For example, the best projection with three 3-level factors is columns $1u^2$, 2 and 5, where $1u^2$ means permutation $u^2$ applies to column 1 in Table 2. It should be mentioned that, with exception of one trivial case in which the full factorial design is the only nonisomorphic design, none of the designs in Tables 4 and 5 is combinatorially isomorphic to the minimum aberration designs given by Xu and Wu (2001).

TABLE 5
*Minimum aberration projections of $L_{18}$ with the two-level factors*

| # of factors | Columns | $(\beta_3, \beta_4, \beta_5)$ | Resolution |
|---|---|---|---|
| 3 | 0 1 2 | (0, 0, 0) | |
| 4 | 0 $1u^2$ 2 5 | (0, 0.5, 1) | IV |
| 5 | 0 1 $3u$ 4 7 | (0, 3.75, 0) | IV |
| 6 | 0 $1u$ 2 $3u^2$ $4u^2$ 7 | (0, 10.0625, 0) | IV |
| 7 | 0 $1u$ 2 $3u^2$ $4u^2$ 5 7 | (1.25, 14.21874, 7.40625) | III |
| 8 | 0 1 $2u$ $3u$ 4 $5u$ 6 7 | (2.5, 22.5, 17.3125) | III |

**5. Concluding remarks and some discussion.** This paper proposes a geometric approach in studying factorial designs with quantitative factors. When factors are quantitative, the traditional mathematical treatment, which is appropriate for designs with nominal factors, no longer applies. The key in our approach is indicator functions and their polynomial forms as expanded to OCBs. They are used to distinguish designs' geometric structures, which carry the designs' properties. This approach is still appropriate even for ordinal categorical factors.

In Corollary 2.2, the connection between $b_{\mathbf{t}}$ and the aliasing of contrasts with no common factors is given. For contrasts that are not disjoint, the calculation of their aliasing (correlation) is more complex, but still depends on the $b_{\mathbf{t}}$'s. In the following, a general formula is offered, which covers the situations of disjoint and nondisjoint contrasts. Let $\{C_{\mathbf{t}}(\mathbf{x}), \mathbf{t} \in \mathcal{T}\}$ be an SOCB. For each factor $X_i$, any product of its two contrasts $C_u^i(x)C_v^i(x)$ can be expressed as a linear combination of $C_0^i(x), C_1^i(x), \ldots$, and $C_{s_i-1}^i(x)$ on the space $\{0, \ldots, s_i - 1\}$. Let

$$C_u^i(x)C_v^i(x) = \sum_{w=0}^{s_i-1} h_w^{(i,u,v)} C_w^i(x) \qquad \text{for } x = 0, 1, \ldots, s_i - 1.$$

The correlation of two contrasts $C_{\mathbf{u}}(\mathbf{x})$ and $C_{\mathbf{v}}(\mathbf{x})$ can be written as a linear combination of $b_{\mathbf{t}}$'s by the following formula:

$$\begin{aligned}
\frac{1}{N} \sum_{\mathbf{x} \in \mathcal{A}} C_{\mathbf{u}}(\mathbf{x}) C_{\mathbf{v}}(\mathbf{x}) &= \frac{1}{N} \sum_{\mathbf{x} \in \mathcal{A}} \prod_{i=1}^{k} C_{u_i}^i(x_i) C_{v_i}^i(x_i) \\
&= \frac{1}{N} \sum_{\mathbf{x} \in \mathcal{A}} \prod_{i=1}^{k} \sum_{w_i=0}^{s_i-1} h_{w_i}^{(i,u_i,v_i)} C_{w_i}^i(x_i) \\
&= \sum_{w_1=0}^{s_1-1} \sum_{w_2=0}^{s_2-1} \cdots \sum_{w_k=0}^{s_k-1} \left( \prod_{i=1}^{k} h_{w_i}^{(i,u_i,v_i)} \right) \left( \frac{1}{N} \sum_{\mathbf{x} \in \mathcal{A}} C_{\mathbf{w}}(\mathbf{x}) \right) \\
&= \sum_{\mathbf{w} \in \mathcal{T}} \left( \prod_{i=1}^{k} h_{w_i}^{(i,u_i,v_i)} \right) b_{\mathbf{w}}.
\end{aligned}$$ (5.1)

Note that for disjoint contrasts $C_{\mathbf{u}}$ and $C_{\mathbf{v}}$, where $\mathbf{u} = (u_1, \ldots, u_k)$ and $\mathbf{v} = (v_1, \ldots, v_k)$,

$$h_{w_i}^{(i,u_i,v_i)} = \begin{cases} 1, & \text{if } w_i = u_i + v_i, \\ 0, & \text{otherwise,} \end{cases}$$ (5.2)

for $i = 1, \ldots, k$, and Corollary 2.2 can be also derived from equations (5.1) and (5.2). To demonstrate the calculation of aliasing between nondisjoint contrasts using (5.1), let us consider the two designs in Example 2.1. Notice that for $x \in \{0,1,2\}$,

$$C_1(x)C_1(x) = \tfrac{1}{\sqrt{2}}C_2(x) + 1,$$
$$C_1(x)C_2(x) = \tfrac{1}{\sqrt{2}}C_1(x) \quad \text{and}$$
$$C_2(x)C_2(x) = \tfrac{-1}{\sqrt{2}}C_2(x) + 1.$$

Therefore, by (5.1) the correlation between two nondisjoint contrasts, say $C_{110}(\mathbf{x})$ and $C_{101}(\mathbf{x})$, equals $(\tfrac{1}{\sqrt{2}}b_{211}+b_{011})/b_{000}$. Its value is $\tfrac{1}{4}$ for (2.4) and $\tfrac{1}{2}$ for (2.5). In general, the aliasing among nondisjoint contrasts has also been captured in the wordlength patterns, which can be viewed as a summary measure of aliasing and are easy to compute. In theory, one can derive a criterion that explicitly calculates aliasing among all pairs of contrasts by laborious computation and this deserves some further investigation.

The $\beta$ wordlength pattern can be generalized when effect order is defined in other ways. For example, the effect order in Definition 4.1 is based on the number of factors that correspond to an effect, whereas the effect order in Definition 4.2 is based on the degree of the polynomial that represents an effect. The two characteristics can be combined to rank effect orders as follows: (a) first use the degree of polynomials to rank effects and then for those effects with the same order, use the number of factors to further rank their order; or (b) first use the number of factors to rank effects and then for those effects with the same order, use the degree of polynomials to further rank their order. For three-level designs, (a) generates the following order:

$$l >> q >> ll >> lq == ql$$
$$>> lll >> qq >> llq == lql == qll >> llll >> \cdots,$$

and for (b), the effect order follows:

$$l >> q >> ll >> lq == ql$$
$$>> qq >> lll >> llq == lql == qll >> llll >> \cdots.$$

In either case, the wordlength patterns can be defined by taking the sum of $(b_\mathbf{t}/b_\mathbf{0})^2$ over coefficients of the $C_\mathbf{t}$'s that are considered equally important. The corresponding aberration criteria sequentially minimize these sums starting from the most important effects. In general, this methodology is very flexible and can be applied on any reasonable effect orders. Note that Corollary 4.1 and Theorem 4.1 still hold under these wordlength patterns. In addition, although $\beta$ wordlength patterns are defined on the coefficients with respect to an OPB, they can be defined with respect to an SOCB as long as an appropriate effect ordering exists.

We chose to present this work in a self-contained fashion rather than with full algebraic geometry language, so that the ideas are more accessible to the statistical community. Nonetheless, this work is another example of how algebraic geometry can be applied to statistics, and we will continue to explore the connections between the two fields. For other applications of algebraic geometry methods in statistics, see Pistone, Riccomagno and Wynn (2000).

## APPENDIX

The geometrically nonisomorphic projected designs of $L_{18}$ are listed in Table 6. All but two pairs of these designs have distinct wordlength patterns as defined in (4.3).

TABLE 6
*Combinatorially and geometrically nonisomorphic projected designs of $L_{18}$*

| Comb. nonisomorphic | Geom. nonisomorphic | WLP $(\beta_3, \beta_4, \beta_5)$ |
|---|---|---|
| 18-3.1 | 1 2 3 | (0.09375, 0.09375, 0.2813) |
|  | $1u^2$ 2 3 | (0, 0.375, 0) |
| 18-3.2 | 1 2 5 | (0.09375, 0.594, 0.281) |
|  | $1u^2$ 2 5 | (0, 0.125, 0.75) |
|  | $1u\ 2u^2$ 5 | (0.375, 0.125, 0.375) |
|  | $1u^2\ 2u^2$ 5 | (0, 0.5, 0) |
| 18-3.3 | 1 3 4 | (0.375, 0.375, 1.125) |
|  | $1u$ 3 4 | (0, 1.5, 0) |
| 18-4.1 | 2 3 4 5 | (0.375, 0.515, 1.313) |
|  | $2u$ 3 4 5 | (0.1875, 0.938, 0.938) |
|  | $2u^2$ 3 4 5 | (0.281, 0.797, 1.406) |
|  | $2u^2\ 3u^2$ 4 5 | (0, 2.064, 0) |
| 18-4.2 | 1 2 3 6 | (0.1875, 0.75, 1.875) |
|  | $1u$ 2 3 6 | (0.375, 0.891, 1.313) |
|  | $1u^2$ 2 3 6 | (0.281, 1.172, 1.031) |
|  | $1u\ 2u^2$ 3 6 | (0.5625, 0.75, 1.125) |
|  | $1\ 2u^2$ 3 6 | (0, 1.875, 0) |
|  | $1\ 2u^2\ 3u$ 6 | (0.281, 0.844, 1.406) |
|  | $1u\ 2u^2\ 3u$ 6 | (0.469, 0.985, 0.844) |
|  | $1u^2\ 2u^2\ 3u$ 6 | (0.656, 0.422, 1.406) |
|  | $1u\ 2u\ 3u$ 6 | (0.1875, 1.5, 0.75) |
|  | $1\ 2u\ 3u$ 6 | (0.1875, 0.9375, 1.313) |
| 18-4.3 | 1 2 3 4 | (0.5625, 0.9375, 1.688) |
|  | $1u$ 2 3 4 | (0.281, 1.781, 0.844) |
|  | $1u\ 2u$ 3 4 | (0, 2.625, 0) |
| 18-4.4 | 1 2 5 6 | (0.5625, 0.9375, 1.688) |
|  | $1u$ 2 5 6 | (0.281, 1.781, 0.844) |
|  | $1u\ 2u^2$ 5 6 | (0.75, 1.125, 0.75) |
|  | $1u\ 2u\ 5u$ 6 | (0.1875, 1.6875, 1.3125) |

**Acknowledgments.** The authors are grateful to the referees for helpful comments.

Institute of Statistical Science
Academia Sinica
Taipei 115
Taiwan

Department of Epidemiology
and Population Health
Albert Einstein College of Medicine
Bronx, New York 10461
USA
e-mail: kyeblue@yahoo.com